\documentclass[11pt]{article}
\usepackage{amsfonts}
\def\square{{\vcenter{\hrule height.4pt
      \hbox{\vrule width.4pt height5pt \hskip5pt
           \vrule width.4pt}
      \hrule height.4pt}}}
\input{epsf}
\def\qed{\hfill$\square$}

\def\aa{\alpha}
\def\bb{\beta}
\def\cc{\gamma}

\def\tt{\tau}

\def\R{{\mathbb R}}

\newtheorem{proposition}{Proposition}[section]
\newtheorem{theorem}[proposition]{Theorem}

\newtheorem{lemma}[proposition]{Lemma}

\title{Untwisting Heegaard diagrams in 3-space}
\author{David Gillman \and Dale Rolfsen}

\begin{document}

\maketitle

\abstract{We show that if $V^3$ is a handlebody in $\R^3$, with curves
$J_1, \dots , J_g \subset \partial V$ which are the attaching curves for a
Heegaard
splitting of a homology sphere, then there exists a homeomorphism
$h\colon V \to V$ so that each of the curves $h(J_i)$ bounds an
orientable surface in $\R^3 - int(V)$.  This leads to a new characterization 
of homology spheres and also
contradicts a remark of Haken \cite{Haken} in 1969 regarding the Poincar\'e
homology sphere.

\section{Definitions and the main theorem.}
All spaces and maps are taken to be in the piecewise-linear category and
we will use standard PL
terminology as, for example, in \cite{Hempel} and \cite{RourkeSanderson}.
Every closed orientable
3-dimensional manifold $M^3$ has a Heegaard splitting $M = V \cup W$ as
the union of two
3-dimensional handlebodies $V$ and $W$, of some genus $g$, which are
disjoint, except that they
have common boundary: $V \cap W = \partial V = \partial W$.  
To specify $M$ it is enough to specify a set of
``Heegaard'' curves $J_1, \dots , J_g$ in the boundary of the handlebody
$V$ to which the meridian
curves of $W$ (which bound disjoint disks in $W$) are attached.  Defining
$J = J_1 \cup \cdots \cup J_g$, the pair $(V,J)$ uniquely determines $M$ 
and will be called a {\it Heegaard diagram} defining $M$.  Indeed, $M$ is
just the result of attaching 2-handles to $V$ along the $J_i$, and then a final
3-handle.  

{\bf Definition:} If $J \subset V \subset
R^3$, we say that a Heegaard diagram $(V,J)$ is {\it untwisted} if each curve of $J$
satisfies: (*) it bounds an orientable surface in $R^3 - int(V)$.  

If the curves $J_i$ bound {\it disjoint} surfaces in $R^3 - int(V)$, we shall
say $(V,J)$ is {\it strongly} untwisted.
\medskip
We now state our main result.

\begin{theorem} \label{main}
Suppose $V^3$ is a handlebody of genus $g$ in $\R^3$ with a system of
curves
$J \subset \partial V$ such that $(V,J)$ is a Heegaard diagram defining a
homology sphere.  Then there is a homeomorphism $h\colon V
\to V$ such that the
diagram $(V,h(J))$ is untwisted.  In particular every Heegaard splitting
of every
homology 3-sphere can be realized by a handlebody embedded in 3-space so
that the Heegaard curves are untwisted. \end{theorem}

Recall that a {\it homology 3-sphere} is a closed 3-manifold $M$ with $H_1(M) = 0$ 
(integer coefficients are assumed for all our homology groups).  If, moreover,
$\pi_1(M) = 1$ the manifold is a {\it homotopy} sphere.
 
Theorem \ref{main} applies to {\it any} embedding of $V$ in
$R^3$, no matter how
horribly knotted and entangled the handles may be.  The surfaces of (*)
may well intersect each other.  Indeed they
cannot be pairwise disjoint, unless the
homology sphere is actually a homotopy sphere
because of the following theorem, due to W. Haken \cite{Haken}.  Another
proof, by E. Rego and C. Rourke, can be found in \cite{RegoRourke}.

\begin{theorem}\label{haken}
{\rm [Haken]}  A closed 3-manifold is a homotopy sphere if and only if it has 
a Heegaard diagram which is strongly untwisted in $\R^3$. \qed
\end{theorem}

Our methods give a similar characterization of homology spheres.

\begin{theorem}\label{char}  A closed 3-manifold is a homology sphere if and only if it 
has a Heegaard diagram which is untwisted in $\R^3$. 
\end{theorem}

We remark that
Theorem \ref{main} contradicts a statement made by Haken in 
\cite{Haken} page 152, where it is claimed that one can {\it not} untwist a
genus 2 Heegaard diagram of the
Poincar\'e dodecahedral homology sphere.  The algorithm used in
\cite{Haken}, pp. 151-2, to attack
this problem is correct.  However, the input variables must be carefully
inserted or the solution
is lost.  In Section \ref{hempel}, we consider precisely this example,
using the Heegaard diagram
of the dodecahedral homology sphere given by Hempel \cite{Hempel} p. 19, which
we reproduce (with permission) as Figure 1.  
The diagram is
certainly {\it not} untwisted ($J_2$ links both of the handles of $V$) but we use
a simple argument to give an explicit untwisted reimbedding (see Figure 4 ahead) of
the Heegaard diagram, as guaranteed by Theorem \ref{main}. but contrary to Haken's
assertion that this is impossible.

\begin{figure}[ht]
{\centering
\vbox{
\epsfysize=2in
\epsffile{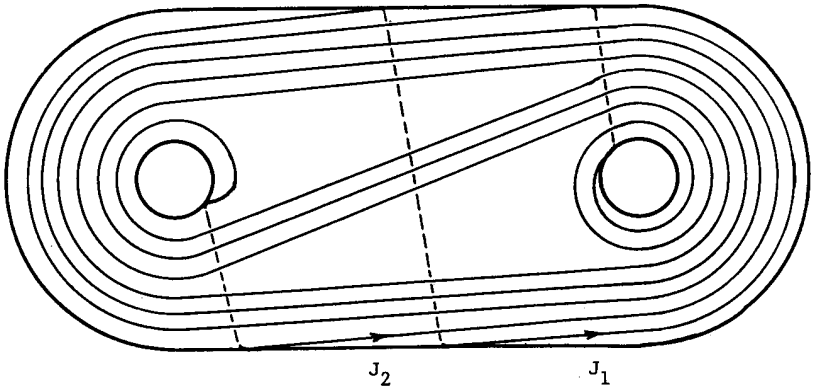} 
}
}
\caption{A Heegaard splitting of Poincar\'e's dodecahedral space.} 
\label{fig1} 
\end{figure} 

\section{Proof of Theorems \ref{main} and \ref{char}}

In this section we assume given a handlebody $V$ embedded in 3-space, and
a set of Heegaard curves $J = J_1 \cup \cdots \cup J_g$ as in the hypothesis
of Theorem \ref{main}.  It is 
convenient to view $V$ as a regular neighbourhood of a 1-dimensional spine 
$X \subset V \subset \R^3$.  We can moreover consider $X$ as a bouquet of $g$
simple closed curves, $X = X_1 \cup \cdots \cup X_g$ which are disjoint except
for their common basepoint $x_0$.  Taking regular neighbourhoods of $x_0$, and
then the $X_i$, defines $V$ as a union of a $0$-handle with $g$ $1$-handles
disjointly attached.  With this convention, we can define homology classes
$\aa_i$ and $\bb_i$, $i = 1,\dots, g$, in $H_1(\partial V)$ by taking $\aa_i$
parallel (in $V$) to the $i$-th handle $X_i$ and $\bb_i$ transverse to that handle,  
with the properties:

(i) $\aa_i$ is represented by a curve on $\partial V$ which is 
parallel in $V$ to $X_i$, in that $\aa_i$ and $X_i$ cobound an annulus in $V$, moreover
(by adding a multiple of $\bb_i$ to $\aa_i$ if necessary)
we may assume $lk(\aa_i,X_i) = 0$, where $lk$ denotes the linking pairing between disjoint
$1$-cycles in $\R^3$;

(ii) $\bb_i$ is represented by the boundary of a (meridian) disk 
in $V$ which is transverse to $X_i$ at a single point and disjoint from 
$X_j, j\ne i$, orientations chosen so that $lk(\bb_i,X_i) = 1$;

(iii) the $\aa_i$ and $\bb_i$ form a symplectic basis for $H_1(\partial V)$,
meaning that the intersection pairing on $\partial V$ satisfies 
$\aa_i \cdot \bb_i = 1 = -\bb_i \cdot \aa_i$ and all other pairs have intersection zero.

Each Heegaard curve $J_i$, if given an orientation,
determines a homology class in $H_1(\partial V)$, expressible uniquely as 
$$J_i = \sum_{j=1}^g (a_{ij} \aa_j + b_{ij} \bb_j),$$
thus defining two square $g$-by-$g$ matrices $A = (a_{ij})$ and $B = (b_{ij})$.
By (i), $H_1(V)$ can be taken to have basis $\aa_1, \dots, \aa_g$, and moreover $A$ can be
considered a presentation matrix for $H_1(M)$.  Because we are assuming that $M$ is a
homology sphere, we must have $\det(A) = \pm 1$, and we may take it to be $+1$ by changing
orientation of one of the $J_i$ if necessary.

Now the description of $V$ as a handlebody has some ambiguities.  Think of a small regular
neighbourhood $V_{\epsilon}$ of the spine $X$.  The handle corresponding to $X_i$ may be
slid over $X_j$ to form a new bouquet $X'$, as illustrated in Figure 2.
Letting
$V'_{\epsilon}$ denote a small regular neighbourhood of $X'$, then the transformation
$V_{\epsilon} \to V'_{\epsilon}$ can be realized by an ambient isotopy within $V$, fixed
on $\partial V$ (see \cite{Rolfsen}, p. 95, for an illustration of such an isotopy).  Since $V$
and $V_{\epsilon}$ cobound a product (surface times interval) the same is true of $V$ and
$V'_{\epsilon}$ and therefore $V$ has $V'_{\epsilon}$ and $X'$ as spines as well.  However, this
changes our choice of basis elements, in that  $\aa_i$ is replaced by $\aa_i \pm \aa_j$.  In
other words, a handle slide multiplies the presentation matrix $A$ by an elementary matrix,
which is by definition the identity matrix but with $\pm 1$ in a single off-diagonal entry.  It
is well-known that any unimodular integer matrix (and its inverse) may be written as a product
of elementary matrices, and so we may assume without loss of generality that $A$ is the identity
matrix.

\begin{figure}[ht]
{\centering
\vbox{
\epsfysize=1.3in
\epsffile{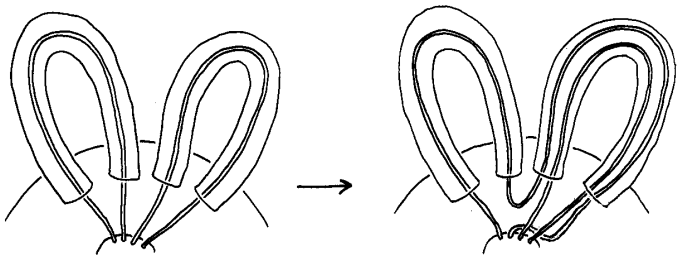} 
}
}
\caption{ Sliding one handle over another.} 
\label{fig2} 
\end{figure}

Therefore, we now assume 

(iv)  $J_i = \aa_i + \sum_{j=1}^g b_{ij} \bb_j.$

It may happen that $\aa_i$ and $X_j$ have nonzero linking number.  However, by a trick
similar to a handle slide, one can swivel the foot of the $i$-th handle around that of the
$j$-th handle appropriately, to ensure the following condition as well:

(v) $lk(\aa_i,X_j) = 0$ for all $i,j \in \{1, \dots, g\}$

Next we look at condition (*), which is well-known to be equivalent to the requirement that
the curve represent a trivial homology class in $R^3 - int(V)$.

\begin{lemma}\label{linking}  The Heegaard diagram $(V,J)$ is untwisted if and only if all
the following linking numbers vanish: $lk(J_i,X_j) = 0$ for all $i,j \in \{1,\dots,g\}$ and
$lk(J_i,J^+_i) = 0$, for all $i \in \{1,\dots,g\}$ where $J^+_i$ is the pushoff of $J_i$
in the outward normal direction away from $\partial V$.
\end{lemma}

\noindent {\bf Proof:}  If the diagram is untwisted, then the existence of the
surfaces guaranteed by (*) shows all the linking numbers are zero.  The converse is also a
standard argument.  For fixed $i$ consider the pushoff curve $J^+_i$, which, as a
knot in $R^3$, bounds a Seifert surface $F^2$ in $\R^3$.  We can assume $F$ is
disjoint from $J_i$ and also that $F$ is in general position relative to the spine
$X$.  By the latter assumption, we can suppose, up to isotopy, that $F$ meets $V$ in
only meridian disks on the 1-handles.  By the assumption that $lk(J_i,X_j) = 0$, the
intersection of $F$ with each handle of $V$ occurs in meridian disks which can
be removed in pairs and replaced by tubes parallel to $X_j$, resulting in a
new surface disjoint from $V$.  Therefore $J^+_i$ is
homologically trivial in $R^3 - int(V)$ and consequently so is each $J_i$; the
diagram is untwisted. \qed

Now we note that the linking numbers can be read from the matrix $B$.  We calculate, with
the help of (i), (ii) and (iv), $lk(J_i,X_j) = lk(\aa_i + \sum_k b_{ik} \bb_k, X_j) =
b_{ij} lk(\bb_j,X_j) = b_{ij}.$  Moreover, the condition $lk(J_i,J^+_i) = 0$ is
equivalent (if $J_i = \aa_i$) to $lk(\aa_i,\aa^+_i) = 0$, which is guaranteed by (i).  We
have shown

\begin{lemma}\label{Bzero}  Under the above assumptions, the Heegaard diagram $(V,J)$ is
untwisted if and only if $B$ is the zero matrix.   \qed 
\end{lemma}

To finish the proof of Theorem \ref{main} it remains to demonstrate how to build a
homeomorphism $h:V \to V$, so that the curves $h(J)$ have matrix $B = 0$, or in other
words, $h(J_i) = \aa_i$ in $H_1(\partial V)$ for all $i$.  The homeomorphism $h$
will be realized as a sequence of Dehn twists along disks properly embedded in $V$.  Recall
that a Dehn twist of a thickened disk, parametrized as $D^2 \times [0,1]$, is the
self-homeomorphism $(z,t) \to (e^{2\pi ti}z,t)$, where $D^2$ is the unit disk in the
complex plane.

We observe that (with our assumptions in force), $B$ is a symmetric matrix.  Indeed
$b_{ij} = lk(J_i,X_j) = lk(J_i,J_j)$, the latter equality due to the fact that $J_j$ is
homologous to $X_j$ in $V$, and so they cobound a 2-chain which may be taken to be
disjoint from $J_i$.  Since $lk(J_i,J_j)) = lk(J_j,J_i)$ we conclude $b_{ij} = b_{ji}$.

Now, suppose an off-diagonal term $b_{ij}$ is not zero.  Consider a
simple closed curve $\tt_{ij}$ which is the connected sum, on $\partial V$ of $\bb_i$ and
$\bb_j$ as shown in Figure 3.  Observe that $\tt_{ij}$ bounds a disk properly embedded in
$V$, the boundary-connected sum of two meridian disks in $V$, with interior pushed into
$int(V)$. Let $T_{ij}: V \to V$ denote the Dehn twist corresponding to this disk in the
appropriate direction so that the effect of $T_{ij}$ on $H_1(\partial V))$ is 

 $$T_{ij}(\cc) = \cc + (\cc \cdot \tt_{ij}) \tt_{ij} = \cc + (\cc
\cdot \bb_i +\cc \cdot \bb_j)(\bb_i + \bb_j).$$

\begin{figure}[ht]
{\centering
\vbox{
\epsfysize=3in
\epsffile{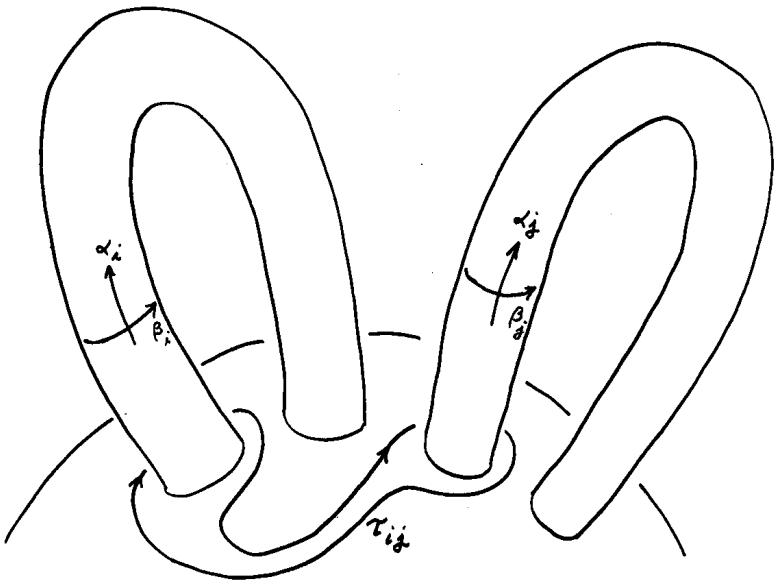} 
}
}
\caption{The curve $\tt_{i,j}$ on the boundary of the handlebody $V$.} 
\label{fig3} 
\end{figure} 
 
In particular, $T_{ij}(J_i) = J_i + \bb_i + \bb_j$ and 
$T_{ij}(J_j) = J_j + \bb_i + \bb_j$, while $T_{ij}(J_k) = J_k$ for $k \ne i, j.$ 
Applying this to the linking numbers, we see that the result of applying $T_{ij}$ is to
increase the entries $b_{ij}$ and $b_{ji}$ by one.  The diagonal entries $b_{ii}$ and
$b_{jj}$ may also change -- we don't care for the moment -- but all the other entries of
$B$ are unchanged by this Dehn twist.  Therefore, by applying an appropriate power
(possibly negative) of $T_{ij}$ to $V$ we obtain a new system of Heegaard curves with
$b_{ij} = b_{ji} = 0$.  We iterate this for all $i,j$ and eventually have a sequence of
Dehn twists which takes the curves $J$ to a set of curves with $B$ matrix diagonal.

Finally, the diagonal entries of $B$ can be modified at will by further Dehn twists
along the meridian disks bounded by the $\bb_i$.  The composite of all these Dehn twists
provides the homeomorphism $h:V \to V$ with $h(J_i) = \aa_i$ in $H_i(\partial V)$,
completing the proof of Theorem \ref{main}. \qed

It remains to prove Theorem \ref{char}, which will follow directly from Theorem \ref{main}
and the following. 

\begin{lemma}
If $M$ has a Heegaard diagram $(V,J)$ which is untwisted in $\R^3$, then $M$ is a
homology sphere.
\end{lemma}

\noindent {\bf Proof:}
The homology $H_1(M)$ can be computed as $H_1(V) / i_*\langle J \rangle$, where
$\langle J \rangle$ is the subgroup of $H_1(\partial V)$ generated by the
homology classes of the curves of $J$, and 
$i_* \colon H_1(\partial V) \to H_1(V)$ is
induced by inclusion.  

The assumption that $(V,J)$ is untwisted in $\R^3$ is equivalent to saying
that $\langle J \rangle$ is a subgroup of the kernel 
$$K := \ker \{H_1(\partial V) \to H_1(\R^3 - int(V))\}.$$
On the other hand, a Mayer-Vietoris argument shows that 
$0 = H_1(\R^3) = H_1(V) / i_*K$, or in other words, $i_*(K) = H_1(V)$.
We now have that $\langle J \rangle$ is a subgroup of $K$ and both are
free abelian of rank $g$.  This alone would not imply  $\langle J
\rangle = K$, but we observe that for any Heegaard splitting whatsoever, 
$\langle J \rangle$ has a complement in $H_1(\partial V)$.  That is,
there is a free abelian subgroup $G$ of $H_1(\partial V)$ so that
$H_1(\partial V) = \langle J \rangle \oplus G$, and both 
$\langle J \rangle$ and $G$ have rank $g$.  (This can be seen, for example,
from the fact that $\partial V$ cut along $J$ becomes a planar surface.)  If $K$
intersected $G$ nontrivially, we would have 
$K \cong \langle J \rangle \oplus K \cap G$,
and obtain a contradiction by a rank argument.  Therefore $\langle J
\rangle = K$ and $i_* \langle J \rangle = i_*K = H_1(\partial V)$, which
establishes that $H_1(M)$ is trivial. \qed

\section{An example: dodecahedral space}\label{hempel}

We now turn to a concrete untwisting of a specific example: the Poincar\'e dodecahedral
homology sphere.  Its fundamental group is finite, of order 120, and perfect
(its abelianization being the trivial group).  A specific genus two Heegaard diagram
$(V,J)$ is given in Figure 1.  Here $V$ is a standardly-embedded genus 2 handlebody in
$\R^3$.   In the terminology of the previous section, the spine $X = X_1 \cup X_2$
can be pictured in the plane of the page, with $X_1$ on the left, $X_2$ on the
right.  The linking number of an oriented curve on $\partial V$ with $X_1$ can be
thought of as the algebraic number of times the curve through the hole ``from front
to back'' on the left side, and likewise for linking with $X_2$ and the right hole.  

Our homeomorphism $h:V \to V$ will be the composite of three Dehn twist
homeomorphisms, using horizontal disks transverse to each of the three pillars of the
diagram.  Reading from left to right, these disks will be twisted $x,y$ and $z$
times, respectively, with the convention that a positive twist corresponds to a
right-hand screw sense.  After these twists are performed, we wish each of the two
Heegaard curves to go through each of the two holes of $V_3$ algebraically zero
times, yielding:

$0 - x - 2y  =  0$  ($J_1$ goes through the left hole zero times)

$0 + 2y + z =  0$   ($J_1$ goes through the right hole zero times)

$-1 + x + y  =  0$  ($J_2$ goes through the left hole zero times)

$-1 - y + 0z =  0$  ($J_2$ goes through the right hole zero times)

\noindent with unique solution: $x=2, y=-1, z=2.$  The corresponding Heegaard diagram
newly embedded in $\R^3$, is pictured in Figure 4 (for clarity, the 
Heegaard surface $\partial V$ is transparent).  

\begin{figure}[ht]
{\centering
\vbox{
\epsfysize=3in
\epsffile{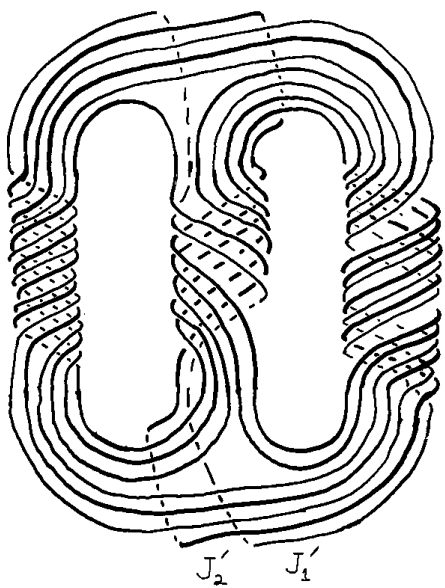} 
}
}
\caption{ An untwisted Heegaard diagram of Poincar\'e's homology sphere.} 
\label{fig4} 
\end{figure}

Notice that in the Heegaard diagram of Figure 1, both curves $J_i$ are, 
considered
individually, unknotted, though they are interlinked.  On the other hand, 
in the untwisted diagram of the same manifold, depicted in 
Figure 4, each of $J'_1 = h(J_1)$ and $J'_2 = h(J_2)$ is a
nontrivial knot.   The curve $J'_1$ is in fact a pretzel knot of type 5,-3,5 
(see Figure 5) with
Alexander  polynomial $1 - 3t + t^2$  and Jones polynomial
$t^{-1} + t^{-3} - t^{-5} + t^{-6} -2t^{-7} + t^{-8} - t^{-9} + t^{-10}$.
Remarkably, $J'_2$ is also a 5, -3, 5 pretzel knot, though differently
embedded in $\partial V$.  

\begin{figure}[ht]
{\centering
\vbox{
\epsfysize=1.5in
\epsffile{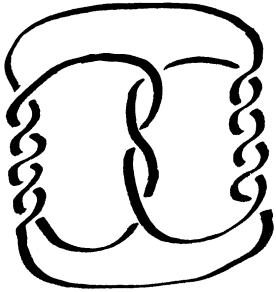} 
}
}
\caption{The 5, -3, 5 pretzel knot.} 
\label{fig5} 
\end{figure} 

\section{Questions and acknowledgements.}

{\bf Question 1:} Why are the two Heegaard curves as pictured in Figure 4
of the same knot type?  Is this typical of untwisted Heegaard diagrams in
$\R^3$?

{\bf Question 2:} In Figure 1, the curve $J_1$ bounds a disk in 
$\R^3 - int(V)$, and $V$ plus a 2-handle attached along $J_1$ is a trefoil
knot exterior (showing, incidentally, why the same manifold is the result of Dehn
surgery along the trefoil).  On
the other hand $J_2$ does not bound in $\R^3 - int(V)$, as it links the
handles.  Is there {\it any} embedding of $V \cup_{J_2} D^2$ into 3-space?

{\bf Question 3:} Must every minimum genus untwisted Heegaard diagram (in $\R^3$) of a
homology sphere ($\ne S^3$) have the attaching curves $J_i$ knotted? 

{\bf Question 4:} Is there a version of Theorem \ref{main} for homotopy spheres?
That is, given a Heegaard diagram $(V,J)$ of a homotopy sphere, and an
arbitrary embedding $V\subset \R^3$, is there a homeomorphism $h \colon V \to V$
so that $(V,h(J))$ is strongly untwisted?  

{\bf Question 5:} Can the study of the intersections of surfaces, associated
with an untwisted Heegaard diagram, lead to invariants of the 3-manifold?
Perhaps the Rohlin or even Casson invariant?
 
We wish to thank John Hempel, Rob Kirby,  Steve Boyer and Zhongmou Li for 
very helpful discussion of aspects of 
this work.  The web utility {\it knotscape} was useful for working out the
example in section \ref{hempel}.

\end{document}